# A non-selfdual 4-dimensional Galois representation

Jasper Scholten

**Abstract**

In this paper it is explained how one can construct non-selfdual 4-dimensional $\ell$-adic Galois representations of Hodge type $h^{3,0} = h^{2,1} = h^{1,2} = h^{0,3} = 1$, assuming a hypothesis concerning the cohomology of a certain threefold. For one such a representation the first 80000 coefficients of its $L$-function are computed, and it is numerically verified that this $L$-function satisfies a functional equation. Also a candidate for the conductor is obtained.

## 1 Introduction

According to the Langlands program there should be a correspondence between cuspidal automorphic representations of $\mathrm{GL}(n)$ and Galois representations in the cohomology of algebraic varieties. Clozel has made this conjecture more precise.

In the case $n = 2$ it is known how to associate a Galois representation to an automorphic representation, due to Eichler, Shimura and Deligne.

For general $n$ Clozel has constructed Galois representations closely related to certain automorphic representations provided that these are selfdual.

For non-selfdual automorphic representations it is not known how to construct corresponding Galois representations. However, there are examples of explicit non-selfdual 3 dimensional automorphic representations, and 3 dimensional Galois representations with the same local $L$-factors as far as computed.

The aim of this paper is to start a search for the same kind of examples for $n = 4$. Assuming a statement about the generators of the $H^2$ of our threefold we find a 4 dimensional Galois representation in the $\ell$-adic $H^3$ of an algebraic threefold. This representation is not isomorphic to a Tate twist of its contragredient representation, even after multiplication with a Dirichlet character. We compute traces of Frobenius at primes of good reduction by counting points on the threefold over finite fields. Moreover, we guess local $L$-factors of the associated $L$-function at primes of bad reduction. This way we can compute the first 80000 coefficients of the $L$-function. Also we guess a conductor of our representation. Having made all these guesses we numerically verify the conjectured funcional equation of the $L$-function. For one particular choice of bad $L$-factors and conductor (namely, the $L$-factors at 2 and 3 are equal to 1, and the conductor is $2^9\,3^9$) we find very convincing evidence that the $L$-function satisfies the functional equation. We also find a value for $L(2)$. It is non-zero. This might be interesting in relation to conjectures of Beilinson and Bloch, which relate the order of the zero of $L(s)$ at $s = 2$ to the rank of the codimension-2 Griffiths group of the threefold.

**Acknowledgements.** I thank Jaap Top for suggesting this topic and helpful conversations, Chad Schoen for his help in some cohomology computations, and Joe Buhler for his programs for evaluating Bessel functions.

## 2  Elliptic surfaces

The constuction in this paper is somewhat similar to the construction of 3-dimensional Galois representations by van Geemen and Top in [10]. They construct representations in the $\ell$-adic $H^2$ of a base changed elliptic surface. We construct representations in the $\ell$-adic $H^3$ of a base changed fibre product of 2 elliptic surfaces. In this section we will review van Geemen and Top's method.

### 2.1  The construction

Let $\pi : \mathcal{E} \to \mathbb{P}^1$ be the Kodaira-Néron model of an elliptic surface, defined over $\mathbb{Q}$. Assume it has a zero section, and non-constant $j$-invariant.

For any variety $V$ defined over a field $k$ we denote the $\ell$-adic cohomology groups $H^i(V_{\bar{k}}, \mathbb{Q}_\ell)$ simply by $H^i(V)$. Define $A(\mathcal{E}) \subset H^2(\mathcal{E})$ to be the subspace generated by a fibre of $\pi$, the image of the 0-section and the fibre components of singular fibres that do not hit the 0-section. The Galois group $G_\mathbb{Q}$ maps $A(\mathcal{E})$ to itself. Denote the Galois module $H^2(\mathcal{E})/A(\mathcal{E})$ by $B(\mathcal{E})$.

The Galois representation on the whole of $H^2(\mathcal{E})$ and also on $B(\mathcal{E})$ will automatically be selfdual because of Poincaré-duality. From now on we choose an embedding $\mathbb{Q}_\ell \to \mathbb{C}$. Our non-selfdual representation will be a $G_\mathbb{Q}$-invariant subspace of $B(\mathcal{E}) \otimes_{\mathbb{Q}_\ell} \mathbb{C}$. The following lemma shows that a sufficient condition for a representation to be non-selfdual is that it has non-real coefficients.

**Lemma 2.1.1** *Let $V/\mathbb{Q}$ be a non-singular variety, and let $M$ be a $d$-dimensional $G_\mathbb{Q}$-invariant subspace of $H^i(V) \otimes_{\mathbb{Q}_\ell} \mathbb{C}$. Let $P(T) = T^d + c_{d-1}T^{d-1} + \ldots + c_0$ be the characteristic polynomial of the $\mathrm{Fr}_q$-action on $M$ for some prime power $q$. Define $\xi$ by $\xi = c_0\,q^{-di/2}$. Then*

$$c_{d-k} = q^{i(k-d/2)}\xi\bar{c}_k.$$

*If $M$ is selfdual then $P(T)$ has real coefficients.*

*Proof.* Write $P(T) = \prod_{j=1}^{d}(T-\alpha_j)$. By Deligne's theorem all $\alpha_j$ have absolute value $q^{i/2}$. So

$$\frac{q^i}{\alpha_j} = \bar{\alpha}_j. \tag{1}$$

Denote by $s_k$ the degree $k$ elementary symmetric polynomial in $d$ variables. The first statement follows from (1) and the identity

$$s_d(T_1, T_2, \cdots, T_d)\, s_k(\frac{1}{T_1}, \frac{1}{T_2}, \ldots, \frac{1}{T_d}) = s_{d-k}(T_1, T_2, \ldots, T_d).$$

If $M$ is selfdual then for each $\mathrm{Fr}_q$-eigenvalue $\alpha$ one has that $q^i/\alpha$ is also an eigenvalue. With (1) this implies that $P(T)$ has real coefficients. □

The way to make a $G_{\mathbb{Q}}$-invariant subspace of $B(\mathcal{E})$ is to let $\mathcal{E}$ have an automorphism of finite order defined over $\mathbb{Q}$. This automorphism will cut the representation on $B(\mathcal{E})$ in several $G_{\mathbb{Q}}$-invariant eigenspaces.

For this define $\pi : \mathcal{E} \to \mathbb{P}^1_s$ to be the Kodaira-Néron model of the base change of an other elliptic surface $\psi : X \to \mathbb{P}^1_t$ over a Galois cover $\phi : \mathbb{P}^1_s \to \mathbb{P}^1_t$ given by

$$t = \phi(s) = \frac{s^3 - 3s^2 + 1}{3s^2 - 3s}.$$

(The sub-index at $\mathbb{P}^1$ is used to denote an affine coordinate). The Galois group of $\phi$ has order 3 and is generated by $\sigma : s \mapsto (s-1)/s$. This map $\sigma$ induces a automorphism on $\mathcal{E}$, which we will also denote by $\sigma$.

From now on $\zeta$ will always denote a fixed primitive third root of unity in $\bar{\mathbb{Q}}_\ell \subset \mathbb{C}$, in $\bar{\mathbb{Q}} = \mathbb{A}^1(\bar{\mathbb{Q}}) \subset \mathbb{P}^1(\bar{\mathbb{Q}})$ or in $\bar{\mathbb{F}}_q = \mathbb{A}^1(\bar{\mathbb{F}}_q) \subset \mathbb{P}^1(\bar{\mathbb{F}}_q)$. From the context it will be clear which $\zeta$ we mean.

The automorphism $\sigma$ induces $G_{\mathbb{Q}}$-equivariant linear automorphisms on the $\ell$-adic cohomology, also denoted by $\sigma$. Since $\sigma$ has order 3, any $\ell$-adic $G_{\mathbb{Q}}$-module $D$ on which $\sigma$ acts decomposes as a sum of three eigenspaces over $\mathbb{Q}_\ell(\zeta)$. Denote the 1-eigenspace of $\sigma|D$ by $D_1$, and the $\zeta$- and $\zeta^2$-eigenspaces of $\sigma|(D \otimes_{\mathbb{Q}_\ell} \mathbb{Q}_\ell(\zeta))$ by $D_\zeta$ and $D_{\zeta^2}$.

## 2.2 Computation of traces of Frobenius

The following lemmas allow us to compute traces of Frobenius acting on an eigenspace of $\sigma$.

**Lemma 2.2.1** *Let $V/\mathbb{F}_q$ be a non-singular variety of dimension $d$ with an automorphism $\sigma$ of order 3, also defined over $\mathbb{F}_q$. Denote the Frobenius automorphism on $V$ and on any $H^j(V)$ by $\mathrm{Fr}_q$. Denote the composition $\mathrm{Fr}_q \circ \sigma^i$ by $\mathrm{Fr}_q^{(i)}$. Let $\mathrm{Fix}(\mathrm{Fr}_q^{(i)}|V)$ denote the number of fixpoints of the map $\mathrm{Fr}_q^{(i)} : V(\bar{\mathbb{F}}_q) \to V(\bar{\mathbb{F}}_q)$, and let $P_j^{(i)}(T) = \det(1 - \mathrm{Fr}_q^{(i)} T|H^j(V))$. Then*

(a)
$$\mathrm{Fix}(\mathrm{Fr}_q^{(i)}|V) = \sum_{j=0}^{2d}(-1)^j \mathrm{Tr}(\mathrm{Fr}_q^{(i)}|H^j(V)).$$

(b) All roots of the $P_j^{(i)}(T)$ have complex absolute value $q^{-j/2}$.

(c) The $P_j^{(i)}(T)$ have coefficients in $\mathbb{Q}$.

(d) $\mathrm{Tr}(\mathrm{Fr}_q^{(i)}|H^0(V)) = 1$ and $\mathrm{Tr}(\mathrm{Fr}_q^{(i)}|H^{2d}(V)) = q^d$.

**Remark** If $i = 0$ then (b) and (c) are part of the Weil conjectures.

*Proof.* (a) is just the Lefschetz trace formula, which we can use because $\mathrm{Fr}_{q^3} = (\mathrm{Fr}_q^{(i)})^3$ has finitely many fixpoints, all of multiplicity 1, hence so has $\mathrm{Fr}_q^{(i)}$.

By Deligne's theorem $(\mathrm{Fr}_q^{(i)})^3|H^j(V)$ has eigenvalues of absolute value $q^{3j/2}$ hence (b) follows.

For (c) we mimic part of the proof of the Weil conjectures as given in [6], chapter 27. First define

$$Z(V, i, T) = \exp\left(\sum_{k=1}^{\infty} \frac{\mathrm{Fix}(\mathrm{Fr}_{q^k}^{(i)}|V)}{k} T^k\right).$$

From (a) it follows that

$$Z(V, i, T) = \frac{P_1^{(i)}(T) \cdots P_{2d-1}^{(i)}(T)}{P_0^{(i)}(T) \cdots P_{2d-1}^{(i)}(T)},$$

so $Z(V, i, T) \in \mathbb{Q}[[T]] \cap \mathbb{Q}_\ell(T)$ hence $Z(V, i, T) \in \mathbb{Q}(T)$ (see [6], 27.9). So we can write

$$Z(V, i, T) = \frac{P(T)}{Q(T)}, \quad P(T) \text{ and } Q(T) \in \mathbb{Q}[T].$$

All $P_j^{(i)}(T)$ have constant coefficient equal to 1, and from (b) it follows that they are relatively prime, hence we may choose $P(T)$ and $Q(T)$ as

$$P(T) = \prod_{j \text{ odd}} P_j^{(i)}(T), \quad Q(T) = \prod_{j \text{ even}} P_j^{(i)}(T).$$

Clearly $Z(V, i, T)$ is independent of $\ell$, hence so are $P(T)$ and $Q(T)$. And for any $j$ the polynomial $P_j^{(i)}(T) \in \mathbb{Q}_\ell[T]$ is characterized as the maximal factor of $P(T)Q(T)$ with roots of absolute value $q^{-j/2}$, so it is also independent of $\ell$. And since it has coefficients in $\mathbb{Q}_\ell$ for every $\ell$, it must have coefficients in $\mathbb{Q}$, and that proves (c).

It is easy to give a direct proof for (d), but it also follows from (c) and the fact that for $j = 0$ and $j = 2d$ we have

$$\mathrm{Tr}(\mathrm{Fr}_q^{(i)}|H^j(V))^3 = \mathrm{Tr}(\mathrm{Fr}_{q^3}|H^j(V)) = q^{3j/2}.$$

□

For the other statements we have to show that $\sigma$ is the identity on $H^0(V)$ and $H^d(V)$. But this is clear since $\sigma$ has order dividing 3, and charateristic polynomial with coefficients in $\mathbb{Q}$, and $\dim H^0(V) = \dim H^d(V) = 1$. □

**Lemma 2.2.2** *Let $D$ be any vectorspace with two commuting endomorphisms $\operatorname{Fr}_q$ and $\sigma$. Suppose $\sigma$ is semisimple and has order dividing 3. Let $D_1$, $D_\zeta$ and $D_{\zeta^2}$ be the three eigenspaces of $\sigma$. Define as before $\operatorname{Fr}_q^{(i)} = \operatorname{Fr}_q \circ \sigma^i$. Then*

$$\operatorname{Tr}(\operatorname{Fr}_q|D_1) = \frac{1}{3}\big(\operatorname{Tr}(\operatorname{Fr}_q|D) + \operatorname{Tr}(\operatorname{Fr}_q^{(1)}|D) + \operatorname{Tr}(\operatorname{Fr}_q^{(2)}|D)\big),$$

$$\operatorname{Tr}(\operatorname{Fr}_q|D_\zeta) = \frac{1}{3}\big(\operatorname{Tr}(\operatorname{Fr}_q|D) + \zeta^2 \operatorname{Tr}(\operatorname{Fr}_q^{(1)}|D) + \zeta \operatorname{Tr}(\operatorname{Fr}_q^{(2)}|D)\big),$$

$$\operatorname{Tr}(\operatorname{Fr}_q|D_{\zeta^2}) = \frac{1}{3}\big(\operatorname{Tr}(\operatorname{Fr}_q|D) + \zeta \operatorname{Tr}(\operatorname{Fr}_q^{(1)}|D) + \zeta^2 \operatorname{Tr}(\operatorname{Fr}_q^{(2)}|D)\big).$$

*Proof.* The lineair map $(I - \zeta\sigma)(I - \zeta^2\sigma) = I + \sigma + \sigma^2$ (with $I$ the identity) is multiplication by 3 on $D_1$ and it maps $D_\zeta$ and $D_{\zeta^2}$ to zero. Hence

$$3\operatorname{Tr}(\operatorname{Fr}_q|D_1) = \operatorname{Tr}(\operatorname{Fr}_q \circ (I + \sigma + \sigma^2)|D)$$
$$= \operatorname{Tr}(\operatorname{Fr}_q|D) + \operatorname{Tr}(\operatorname{Fr}_q^{(1)}|D) + \operatorname{Tr}(\operatorname{Fr}_q^{(2)}|D).$$

The other two equalities follow from similar arguments applied to the maps $(I - \sigma)(I - \zeta\sigma) = I + \zeta^2\sigma + \zeta\sigma^2$ and $(I - \sigma)(I - \zeta^2\sigma) = I + \zeta\sigma + \zeta^2\sigma^2$. □

The previous lemma allows us to compute $\operatorname{Tr}(\operatorname{Fr}_q|B(\mathcal{E})_\zeta)$ if we know the traces $\operatorname{Tr}(\operatorname{Fr}_q^{(i)}|B(\mathcal{E}))$ for $i = 0, 1$ and $2$. The latter traces can be computed using lemma 2.2.1 if we know $\operatorname{Fix}(\operatorname{Fr}_q^{(i)}|\mathcal{E})$ and the $\operatorname{Fr}_q^{(i)}$-action on $A(\mathcal{E})$, since $h^1(\mathcal{E}) = h^3(\mathcal{E}) = 0$. A convenient way of doing this in practise is by means of the following lemma.

**Lemma 2.2.3** *Assume $\mathcal{E}$ has only bad fibres of Kodaira type $I_n$, $II$, $III$ and $IV$. Define the Weierstrass model $\mathcal{E}_{\mathrm{We}} \to \mathbb{P}^1$ to be a complete model for $\mathcal{E} \to \mathbb{P}^1$ given locally by a minimal projective Weierstrass equation (with coefficients in the function field of $\mathbb{P}^1$). Then*

$$\operatorname{Tr}(\operatorname{Fr}_q^{(i)}|B(\mathcal{E})) = \operatorname{Fix}(\operatorname{Fr}_q^{(i)}|\mathcal{E}_{\mathrm{We}}) - (q+1)^2.$$

*Proof.* By lemma 2.2.1 we have

$$\operatorname{Tr}(\operatorname{Fr}_q^{(i)}|B(\mathcal{E})) = \operatorname{Fix}(\operatorname{Fr}_q^{(i)}|\mathcal{E}_{\mathrm{We}}) - (q+1)^2 +$$
$$\big(\operatorname{Fix}(\operatorname{Fr}_q^{(i)}|\mathcal{E}) - \operatorname{Fix}(\operatorname{Fr}_q^{(i)}|\mathcal{E}_{\mathrm{We}})\big) - \big(\operatorname{Tr}(\operatorname{Fr}_q^{(i)}|A(\mathcal{E})) - 2q\big),$$

so we have to show that $\operatorname{Fix}(\operatorname{Fr}_q^{(i)}|\mathcal{E}) - \operatorname{Fix}(\operatorname{Fr}_q^{(i)}|\mathcal{E}_{\mathrm{We}}) = \operatorname{Tr}(\operatorname{Fr}_q^{(i)}|A(\mathcal{E})) - 2q$.

Recall that $A(\mathcal{E})$ is generated by a fibre, the zero section and fibre components that do not hit the zero section. These generators form a basis for $A(\mathcal{E})$. Note that $\sigma$ permutes this basis (since it permutes fibres, and leaves the zero section fixed, and consequently it permutes fibre components not hitting the zero section). Moreover, $\mathrm{Fr}_q$ permutes the basis, and multiplies it with $q$. Hence $\mathrm{Tr}(\mathrm{Fr}_q^{(i)}|A(\mathcal{E}))$ is $q$ times the number of basiselements fixed by $\mathrm{Fr}_q^{(i)}$ settheoretically.

The map $\mathcal{E} \to \mathcal{E}_{\mathrm{We}}$ is a resolution of singularities, and the fibre above each singular point of $\mathcal{E}_{\mathrm{We}}$ is a chain of $\mathbb{P}^1$'s, say with components $F_1, \ldots, F_n$ such that $(F_i, F_{i+1}) = 1$ for $1 \leq i \leq n-1$, and no other intersections. In fact, the fibres in this resolution precisely consist of the components of fibres of the elliptic fibration not intersecting the zero-section. Components that are not fixed by $\mathrm{Fr}_q^{(i)}$ do not contribute to $\mathrm{Fix}(\mathrm{Fr}_q^{(i)}|\mathcal{E}) - \mathrm{Fix}(\mathrm{Fr}_q^{(i)}|\mathcal{E}_{\mathrm{We}})$, and since they are permuted with components (in possibly other fibres) by the $\mathrm{Fr}_q$-action they also don't contribute to $\mathrm{Tr}(\mathrm{Fr}_q^{(i)}|A(\mathcal{E}))$.

Clearly, none of the components in a fibre above a singular point of $\mathcal{E}_{\mathrm{We}}$ that is not fixed by $\mathrm{Fr}_q^{(i)}$ are fixed by $\mathrm{Fr}_q^{(i)}$. The fibre above an singular point fixed by $\mathrm{Fr}_q^{(i)}$ is fixed by $\mathrm{Fr}_q^{(i)}$, and there are only 2 possible $\mathrm{Fr}_q^{(i)}$-actions on the components. It either acts trivial, or it acts as $\mathrm{Fr}_q^{(i)}(F_i) = F_{n+1-i}$. In both cases the part of the fibre consisting of components fixed by $\mathrm{Fr}_q^{(i)}$ is connected. Suppose there are $n'$ such components. Then the fibre contributes $n'q$ to $\mathrm{Tr}(\mathrm{Fr}_q^{(i)}|A(\mathcal{E}))$, and it contributes $n'q + 1$ points to $\mathrm{Fix}(\mathrm{Fr}_q^{(i)}|\mathcal{E})$, instead of the one $\mathrm{Fr}_q^{(i)}$-fixed singular point of $\mathcal{E}_{\mathrm{We}}$. So this tells us that the fibre components contribute equal amount to $\mathrm{Fix}(\mathrm{Fr}_q^{(i)}|\mathcal{E}) - \mathrm{Fix}(\mathrm{Fr}_q^{(i)}|\mathcal{E}_{\mathrm{We}})$ and $\mathrm{Tr}(\mathrm{Fr}_q^{(i)}|A(\mathcal{E}))$. To finish the proof, note that the zero section and a fibre contribute $2q$ to $\mathrm{Tr}(\mathrm{Fr}_q^{(i)}|A(\mathcal{E}))$ and have no contribution to $\mathrm{Fix}(\mathrm{Fr}_q^{(i)}|\mathcal{E}) - \mathrm{Fix}(\mathrm{Fr}_q^{(i)}|\mathcal{E}_{\mathrm{We}})$. $\square$

The previous lemmas allow us to express $\mathrm{Tr}(\mathrm{Fr}_q|B(\mathcal{E})_\zeta)$ in terms of the numbers $\mathrm{Fix}(\mathrm{Fr}_q^{(i)}|\mathcal{E}_{\mathrm{We}})$ for $i = 0, 1$ and $2$. The next lemma shows how to compute these numbers. First we introduce new notation. Let $X_{\mathrm{We},t}$ denote the fibre above $t$ of a Weierstrass model for $X$.

**Lemma 2.2.4** *Define $K^{(i)} \subset \mathbb{P}^1_t - \{\zeta, \zeta^2\}$ by*

$$
\begin{aligned}
K^{(0)} &= \{t \in \mathbb{F}_q \mid x^3 - 3(t+1)x^2 + 3tx + 1 \text{ divides } x^q - x\} \cup \{\infty\}, \\
K^{(2)} &= \{t \in \mathbb{F}_q \mid x^3 - 3(t+1)x^2 + 3tx + 1 \text{ divides } x^{q+1} - x + 1\}, \\
K^{(1)} &= \mathbb{F}_q - K^{(0)} - K^{(1)} - \{\zeta, \zeta^2\}.
\end{aligned}
$$

*Let $R^{(i)}$ denote the total number of points in the fibres of $\mathcal{E}_{\mathrm{We}}$ above $-\zeta$ and $-\zeta^2$ (the ramification points of $\phi$) fixed by $\mathrm{Fr}_q^{(i)}$. Then*

$$\mathrm{Fix}(\mathrm{Fr}_q^{(i)}|\mathcal{E}_{\mathrm{We}}) = 3 \sum_{t \in K^{(i)}} \#X_{\mathrm{We},t}(\mathbb{F}_q) + R^{(i)}.$$

*Proof.* For $t \in \mathbb{P}^1 - \{\zeta, \zeta^2\}$ let $s_1$, $s_2$ and $s_3$ be the points in the preimage $\phi^{-1}(t)$. If $t \neq \infty$ they are precisely the zeroes of the polynomial $x^3 - 3(t+1)x^2 + 3tx + 1$.

If $t \in K^{(0)}$ then the $s_j$ are $\mathbb{F}_q$-rational, and the base change $\mathcal{E}$ has three $\mathbb{F}_q$-rational fibres above the fibre $X_t$, all isomorphic to $X_t$. Hence they contribute $3\#X_t(\mathbb{F}_q)$ to $\mathrm{Fix}(\mathrm{Fr}_q^{(0)}|\mathcal{E})$.

The map $\sigma : s \mapsto (s-1)/s$ permutes the fibre elements $s_1$, $s_2$ and $s_3$. If it acts as Frobenius then $s_j^q = \sigma(s_j)$ hence $s_j^{q+1} - s_j + 1 = 0$, so $s_j \in K^{(2)}$. In this case $\mathrm{Fr}_q^{(2)}$ acts trivial on the three points of $\mathcal{E}(\mathbb{F}_{q^3})$ above each point of $X_t(\mathbb{F}_q)$. And these points contribute to $\mathrm{Fix}(\mathrm{Fr}_q^{(2)}|\mathcal{E})$.

If $t \in K^{(1)}$ then the only possible $\sigma$-action on $\{s_1, s_2, s_3\}$ left is $\sigma(s_j) = s_j^{q^2}$, so points in $\mathcal{E}(\mathbb{F}_{q^3})$ above $X_t(\mathbb{F}_q)$ contribute to $\mathrm{Fix}(\mathrm{Fr}_q^{(1)}|\mathcal{E})$.

If $\zeta \in \mathbb{F}_q$ then the ramification points of $\phi$ are in $\mathbb{F}_q$, and the points in fibre above ramification must also be added. □

## 2.3 Dimension computations

The Betti numbers of elliptic surfaces are computed in [5]. For an elliptic surface with base $\mathbb{P}^1$ we have $h^1 = 0$ and $h^2 = 12(p_a+1)-2$ with $p_a$ the arithmetic genus of the surface. The arithmetic genus can be computed using Tate's algorithm [9] and theorem 12.2 in [5]. The dimension of $A(\mathcal{E})$ can be computed with the Tate algorithm, and $\dim B(\mathcal{E}) = h^2 - \dim A(\mathcal{E})$.

Here we compute the dimensions of the eigenspaces $B(\mathcal{E})_{\zeta^i}$. Since $\sigma$ is defined over $\mathbb{Q}$ the dimensions of $B(\mathcal{E})_\zeta$ and $B(\mathcal{E})_{\zeta^2}$ are equal. So we only need to know $\dim B(\mathcal{E})_1$.

Note that any $x \bmod A(\mathcal{E}) \in B(\mathcal{E})_1$ lifts to the $\sigma$-invariant element $\frac{1}{3}(x + \sigma(x) + \sigma^2(x)) \in H^2(\mathcal{E})_1$. So the canonical injective map

$$H^2(\mathcal{E})_1/A(\mathcal{E})_1 \longrightarrow B(\mathcal{E})_1$$

is also surjective.

The space of invariants $H^2(\mathcal{E})_1$ is via pull-back isomorphic to $H^2(\mathcal{E}/\langle\sigma\rangle)$. If base-change map $\phi$ only ramifies at smooth fibres then $\mathcal{E}/\langle\sigma\rangle$ equals the elliptic surface $X$. In this case $A(\mathcal{E})_1 \cong A(X)$ and $\dim A(\mathcal{E})_1 = \frac{1}{3}(\dim A(\mathcal{E}) - 2) + 2$.

If $\phi$ ramifies at singular fibres then $X$ is the Kodaira-Néron model of $\mathcal{E}/\langle\sigma\rangle$. In order to relate $H^2(\mathcal{E}/\langle\sigma\rangle)$ to the known $H^2(X)$ one should compare the contributions of the singular fibre components of $X$ and of $\mathcal{E}/\langle\sigma\rangle$ to the $H^2$.

For example, if $X$ has fibres of type $IV$ at the ramification locus of $\phi$ then $\mathcal{E}$ has smooth elliptic curves as fibres there, and the quotient $\mathcal{E}/\langle\sigma\rangle$ has irreducible fibres. These irreducible fibres don't contribute to the $H^2$ since they are homologically equivalent to any other fibre. And the two type $IV$-fibres of $X$ both have 3 components and contribute $2(3-1) = 4$ to $\dim H^2(X)$. Hence $\dim H^2(\mathcal{E}/\langle\sigma\rangle) = \dim H^2(X) - 4$. One can reason similarly to find that $\dim A(\mathcal{E})_1 = \dim A(X) - 4$.

## 2.4 Examples

In [10] several examples of non-selfdual 3 dimensional Galois representations in the cohomology of elliptic surfaces are constructed. It turns out that with elliptic surfaces one has pretty good control over the dimensions of several pieces of cohomology. One can construct a non-selfdual representation of any (not too large) dimension, just by playing with the singular fibre configuration. To illustrate this we give examples of such representations of dimensions from 4 to 12.

Let $X \to \mathbb{P}^1_t$ be a rational semistable elliptic surface with $n$ singular fibres $I_{m_1}, \ldots, I_{m_n}$, and smooth fibres above the ramification points of $\phi$. Then $h^2(X) = 10$, $h^2(\mathcal{E}) = 34$, $\sum m_i = 12$ and $\sum(m_i - 1) = 12 - n$. Hence $\dim A(X) = 14 - n$ and $\dim A(\mathcal{E}) = 2 + 3\sum(m_i - 1) = 38 - n$. So $\dim B(X) = n - 4$, $\dim B(\mathcal{E}) = 3n - 4$ and

$$\dim B(\mathcal{E})_\zeta = \frac{1}{2}\left((2n-4) - (n-4)\right) = n.$$

For every $n$, $4 \leq n \leq 12$ we will write down a Weierstrass equation for a surface $X$, and the characteristic polynomial of $\text{Fr}_p$ acting on $B(\mathcal{E})_\zeta$ for some prime $p$ of good reduction. This characteristic polynomial was computed as follows. First the traces $\text{Tr}(\text{Fr}^{(i)}_{p^k}|B(\mathcal{E})_\zeta)$ for $0 \leq i \leq 2$ and $1 \leq k \leq \lceil \frac{n}{2} \rceil$ were computed, as explained above. The coefficients $c_{n-k}$ of the polynomial can be expressed in terms of these traces. The other coefficients were determined with lemma 2.1.1.

In the examples we give the characteristic polynomials are irreducible over $\mathbb{Q}(\zeta)$, so the representations are also irreducible over this field.

**Example** $n = 4$.
Weierstrass equation $y^2 = x^3 - (3t^4 + 24t^3 - 48t + 48)x$
$\qquad\qquad -2t^6 - 24t^5 - 48t^4 + 112t^3 - 48t^2 + 192t - 128$.
Char.pol. at $p = 5$: $T^4 + (1 + 2\zeta)T^3 - 20T^2 - (25 + 50\zeta)T + 625$.

**Example** $n = 5$.
Weierstrass equation $y^2 = x^3 - (219t^4 + 960t^2 + 768)x$
$\qquad\qquad -(1190t^6 + 8304t^4 + 15360t^2 + 8192)$.
Char.pol. at $p = 7$: $T^5 + (8 + 9\zeta)T^4 + (35\zeta - 21)T^3 - (147 + 392\zeta)T^2$
$\qquad\qquad +(2744 - 343\zeta)T + 16807 + 16807\zeta$.

**Example** $n = 6$.
Weierstrass equation $y^2 + txy + y = x^3 + x^2 + x$.
Char.pol. at $p = 2$: $T^6 + (2 + \zeta)T^5 + (4 + 2\zeta)T^4 + (8 + 4\zeta)T^3$
$\qquad\qquad +(16 + 8\zeta)T^2 + (32 + 16\zeta)T + 64\zeta + 64$.

**Example** $n = 7$.
Weierstrass equation $y^2 = x^3 - (16t^4 + 24t^2 + 27)x - (64t^4 + 72t^2 + 54)$.
Char.pol. at $p = 7$: $T^7 - (5 + 5\zeta)T^6 + (7\zeta - 14)T^5 - (245 + 49\zeta)T^4$
$\qquad\qquad +(1715 + 1372\zeta)T^3 + (4802 + 7203\zeta)T^2 + 84035T -$
$\qquad\qquad (823543 + 823543\zeta)$.

**Example** $n = 8$.

Weierstrass equation $y^2 + txy + y = x^3 + t$.
Char.pol at $p=2$: $X^8 - \zeta X^7 + 2 X^6 - 4\zeta X^5 + 16 X^3 - 32\zeta X^2 + 64 X$
$-256\zeta$.

**Example** $n = 9$.
Weierstrass equation $y^2 + txy + (1+t^2)y = x^3 + (1+t-t^2)x + t$.
Char.pol at $p=2$: $T^9 - 3\zeta T^8 - (4+4\zeta)T^7 - 4 X^6 + (8+8\zeta)T^5 - 16\zeta T^4$
$-32 T^3 + 128\zeta T^2 + (384 + 384\zeta)T + 512$.

**Example** $n = 10$.
Weierstrass equation $y^2 + 3txy + (1+t^3)y = x^3 + t$.
Char.pol at $p=2$: $T^{10} + \zeta T^9 + 2 T^8 + 4\zeta T^7 + (8\zeta + 8)T^6 + (32 + 32\zeta)T^4$
$+64 T^3 + 128\zeta T^2 + 256 T + 1024\zeta$.

**Example** $n = 11$.
Weierstrass equation $y^2 + txy + y = x^3 + t^3 x + t$.
Char.pol at $p=2$: $T^{11} - \zeta T^{10} + 2 T^9 - 4\zeta T^8 + 16 T^6 + (-32 - 32\zeta)T^5$
$-128\zeta T^3 + (-256\zeta - 256)T^2 - 512\zeta T - 2048 - 2048\zeta$.

**Example** $n = 12$.
Weierstrass equation $y^2 + txy + (1-7t^3)y = x^3 + (t + t^2 + t^3)x + 1 - 10 t^6$.
Char.pol at $p=2$: $T^{12} + (-2 - \zeta)T^{11} + (4 + 4\zeta)T^{10} + (-8\zeta - 4)T^9$
$+8\zeta T^8 - 32 T^6 - (128 + 128\zeta)T^4 + (256 + 512\zeta)T^3$
$-1024\zeta T^2 + (-1024 + 1024\zeta)T + 4096$.

# 3 Elliptic Threefolds.

Let $\pi : \mathcal{E} \to \mathbb{P}^1$ and $\pi' : \mathcal{E}' \to \mathbb{P}^1$ be two elliptic surfaces such that for every $s$ where one of the surfaces has an additive fibre the other surface has a smooth fibre. Define $W = \mathcal{E} \times_{\mathbb{P}^1} \mathcal{E}'$ and $f : W \to \mathbb{P}^1$. The only singular points on $W$ are nodes in fibres above $s$ where both $\mathcal{E}$ and $\mathcal{E}'$ have multiplicative fibres. There are $nm$ of these singular points in each fibre of Kodaira type $I_n \times I_m$. Let $\tau : \widetilde{W} \to W$ be the blowup at these points and $\tilde{f} = f \circ \tau$. Then $\widetilde{W}$ is nonsingular, and the exceptional divisor of the blowup consists of a $\mathbb{P}^1 \times \mathbb{P}^1$ above each node of $W$.

## 3.1 Dimension computations

In this section we try to determine Betti numbers $h^i$ and dimensions of eigenspaces of the cohomology of $\widetilde{W}$.

First we compute the Euler-Poincaré characteristic. Denote by $m_x$ the number of components in $\pi^{-1}(x)$ if this fibre is singular. Set $m_x = 0$ if $\pi^{-1}(x)$ is non singular. Similarly, define $m'_x$ to be the number of singular fibre components of $\pi'^{-1}(x)$.

**Lemma 3.1.1** *The Euler-Poincaré characteristic $\chi(\widetilde{W})$ of $\widetilde{W}$ is given by*

$$\chi(\widetilde{W}) = 4 \sum_{x \in \mathbb{P}^1} m_x \, m'_x.$$

*Proof.* It is well known that a fibre of Kodaira type $I_n$ has $\chi(I_n) = n$. This includes the case $n = 0$, if $I_0$ stands for a smooth elliptic curve. Hence $\chi(I_n \times I_m) = nm$. This allows us to compute $\chi(W)$ using the fibration $f : W \to \mathbb{P}^1$.

$$\begin{aligned} \chi(W) &= \chi(\mathbb{P}^1)\chi(I_0 \times I_0) + \sum_{x \in \mathbb{P}^1} \left( \chi(f^{-1}(x)) - \chi(I_0 \times I_0) \right) \\ &= \sum_{x \in \mathbb{P}^1} m_x \, m'_x. \end{aligned}$$

The exceptional fibre of the blowup $\tau : \widetilde{W} \to W$ consists of a $\mathbb{P}^1 \times \mathbb{P}^1$ above each node, so

$$\begin{aligned} \chi(\widetilde{W}) &= \chi(W) + \sum_{\substack{w \in W \\ w \text{ node}}} \left( \chi(\mathbb{P}^1 \times \mathbb{P}^1) - 1 \right) \\ &= \chi(W) + \sum_{x \in \mathbb{P}^1} m_x \, m'_x \, (2^2 - 1) = 4 \sum_{x \in \mathbb{P}^1} m_x \, m'_x. \end{aligned}$$

□

**Lemma 3.1.2** $h^{3,0}(\widetilde{W}) = h^{2,0}(\mathcal{E}) + h^{2,0}(\mathcal{E}') + 1.$

*Proof.* We have $h^{3,0}(\widetilde{W}) = \dim H^3(\widetilde{W}, \mathcal{O}_{\widetilde{W}})$. From [4], III, 11.4 and 12.9 it follows that $\tau_* \mathcal{O}_{\widetilde{W}} \cong \mathcal{O}_W$ and $R^i \tau_* \mathcal{O}_{\widetilde{W}} = 0$ for $i > 0$. Hence with the Leray spectral sequence $H^p(W, R^q \tau_* \mathcal{O}_{\widetilde{W}}) \Rightarrow H^{p+q}(\widetilde{W}, \mathcal{O}_{\widetilde{W}})$ it follows that

$$\dim H^3(\widetilde{W}, \mathcal{O}_{\widetilde{W}}) = \dim H^3(W, \mathcal{O}_W)$$

.

The Künneth formula and the fact that $\pi$ and $\pi'$ have 1-dimensional fibres gives

$$R^2 f_* \mathcal{O}_W = R^1 \pi_* \mathcal{O}_\mathcal{E} \otimes R^1 \pi'_* \mathcal{O}_{\mathcal{E}'}.$$

From [4], III, 12.9 and the fact that the fibres of $\pi$ have arithmetic genus 1 it follows that $R^1 \pi_* \mathcal{O}_\mathcal{E}$ is an invertible sheaf on $\mathbb{P}^1$.

From Riemann-Roch it follows that

$$\dim H^0(\mathbb{P}^1, R^1 \pi_* \mathcal{O}_\mathcal{E}) - \dim H^1(\mathbb{P}^1, R^1 \pi_* \mathcal{O}_\mathcal{E}) = \deg(R^1 \pi_* \mathcal{O}_\mathcal{E}) + 1.$$

From $h^1(\mathcal{E}) = 0$ and the Leray spectral sequence

$$H^p(\mathbb{P}^1, R^q \pi_* \mathcal{O}_\mathcal{E}) \;\Rightarrow\; H^{p+q}(\mathcal{E}, \mathcal{O}_\mathcal{E})$$

it follows that $H^0(\mathbb{P}^1, R^1\pi_*\mathcal{O}_\mathcal{E}) = 0$ and
$$\dim H^1(\mathbb{P}^1, R^1\pi_*\mathcal{O}_\mathcal{E}) = \dim H^2(\mathcal{E}, \mathcal{O}_\mathcal{E}) = h^{2,0}(\mathcal{E}).$$

Hence
$$\deg R^1\pi_*\mathcal{O}_\mathcal{E} = -h^{2,0}(\mathcal{E}) - 1.$$

Similarly, one has
$$\deg R^1\pi'_*\mathcal{O}_{\mathcal{E}'} = -h^{2,0}(\mathcal{E}') - 1,$$

and consequently
$$\deg R^2 f_*\mathcal{O}_W = -h^{2,0}(\mathcal{E}) - h^{2,0}(\mathcal{E}') - 2.$$

So $H^0(\mathbb{P}^1, R^2 f_*\mathcal{O}_W) = 0$ and Riemann-Roch implies
$$\dim H^1(\mathbb{P}^1, R^2 f_*\mathcal{O}_W) = h^{2,0}(\mathcal{E}) + h^{2,0}(\mathcal{E}') + 1.$$

Using the Leray spectral sequence $H^p(\mathbb{P}^1, R^q f_*\mathcal{O}_W) \Rightarrow H^{p+q}(\mathbb{P}^1, \mathcal{O}_W)$ it follows that
$$\dim H^3(W, \mathcal{O}_W) = h^{2,0}(\mathcal{E}) + h^{2,0}(\mathcal{E}') + 1.$$

□

**Lemma 3.1.3** *The dimension of $H^1(\widetilde{W})$ is 0.*

*Proof.* In this proof we put the complex topology on our algebraic varieties, and $\mathbb{C}$ denotes the sheaf of locally constant functions. The fibres of $\tau$ are connected, and $H^1(\tau^{-1}(x), \mathbb{C}) = 0$ for all $x$. So $\tau_*\mathbb{C} = \mathbb{C}$ and $R^1\tau_*\mathbb{C} = 0$. The fibres of $\pi$ and $\pi'$ are connected so $\pi_*\mathbb{C} = \mathbb{C}$ and $\pi'_*\mathbb{C} = \mathbb{C}$. Hence by the Künneth theorem we have
$$R^1 f_*\mathbb{C} = R^1\pi_*\mathbb{C} + R^1\pi'_*\mathbb{C}.$$

By [3] one has $H^0(\mathbb{P}^1, R^1\pi_*\mathbb{C}) = H^0(\mathbb{P}^1, R^1\pi'_*\mathbb{C}) = 0$ and consequently
$$H^0(\mathbb{P}^1, R^1 f_*\mathbb{C}) = 0.$$

There is an exact sequence
$$0 \to E_2^{1,0} \to E^1 \to E_2^{0,1} \tag{2}$$

associated to a spectral sequence $E_2^{p,q} \Rightarrow E^{p+q}$. Applying this to the spectral sequence
$$R^p f_*(R^q\tau_*\mathbb{C}) \Rightarrow R^{p+q}\tilde{f}_*\mathbb{C}$$

and taking global sections yields $H^0(\mathbb{P}^1, R^1\tilde{f}_*\mathbb{C}) = 0$.

Clearly $\tilde{f}_*\mathbb{C} = \mathbb{C}$ and $H^1(\mathbb{P}^1, \tilde{f}_*\mathbb{C}) = 0$. So applying (2) to the Leray spectral sequence $H^p(\mathbb{P}^1, R^q\tilde{f}_*\mathbb{C}) \Rightarrow H^{p+q}(\widetilde{W}, \mathbb{C})$ yields $H^1(\widetilde{W}, \mathbb{C}) = 0$. □

Let $r : \widetilde{W} \to \mathcal{E}$ and $r' : \widetilde{W} \to \mathcal{E}'$ be the two canonical projection maps.

**Assumption/Conjecture 3.1.4**

(1) The intersection $r^*(H^2(\mathcal{E})) \cap r'^*(H^2(\mathcal{E}')) \subset H^2(\widetilde{W})$ is one-dimensional and generated by the class of a fibre of $\tilde{f}$.

(2) The space
$$H^2(\widetilde{W})/\bigl(r^*(H^2(\mathcal{E})) + r'^*(H^2(\mathcal{E}'))\bigr)$$
is generated by the following divisors:

– The $\sum_{x \in \mathbb{P}^1} m_x \cdot m'_x$ quadrics that form the exceptional fibre of $\tau$.

– The proper transforms in the blowup $\tau$ of the $\sum_{x \in \mathbb{P}^1}(m_x-1)\cdot(m'_x-1)$ divisors of type $F_{x,i} \times F'_{x,i'}$ with $1 \leq i \leq m_x-1$ and $1 \leq i' \leq m'_x-1$ where $F_{x,i}$ and $F'_{x,i'}$ denote the fibre components of $\pi^{-1}(x)$ and $\pi'^{-1}(x)$ that do not hit the 0-section.

– One other divisor in the case that $\mathcal{E}$ and $\mathcal{E}'$ are isogenous. This means that there is map $\mathcal{E} \to \mathcal{E}'$, commuting with $\pi$ and $\pi'$, and not constant on the fibres. In this case the graph of the isogeny in the fibred product is an extra divisor.

There is no linear relation between these divisors.

**Remark.** Although I don't have a proof for this I believe that it is true. For several threefolds I have computed local L-factors of the Galois representation on the $H^3$ by counting points over many finite fields, and by using the above conjecture to determine the Galois action on $H^2$. All computations showed that 3.1.4 was compatible with the Weil conjectures. Moreover, in the case that $W$ is the fibre product of a semistable elliptic surface with itself the Betti numbers given by the conjecture agree with the Betti numbers of such threefolds as computed by C. Schoen in [7]. In trying to prove the statement I looked at the Leray spectral sequence $H^p(\mathbb{P}^1, R\tilde{f}^q_*\mathbb{Q}) \Rightarrow H^{p+q}(\widetilde{W}, \mathbb{Q})$ which degenerates at $E_2$ (cf.[11]). I can find explicit generators for $H^1(\mathbb{P}^1, R^1\tilde{f}_*\mathbb{Q})$ and $H^2(\mathbb{P}^1, \tilde{f}_*\mathbb{Q})$ but I am not able to compute $H^0(\mathbb{P}^1, R^2\tilde{f}_*\mathbb{Q})$.

We want $\widetilde{W}$ to have an autmorphism $\sigma$ that cuts the Galois representations in pieces. For this we let $\mathcal{E}$ and $\mathcal{E}'$ be Kodaira-Néron models of a cubic basechange of elliptic surfaces $\psi : X \to \mathbb{P}^1_t$ and $\psi' : X' \to \mathbb{P}^1_t$, as in chapter 2.

From now on we will fix the singular fibre types of $W$. The arguments that follow can be done for arbitrary $W$, but they are a bit tedious to write down in full generality, so we restrict ourselves to the case we need. We take $X$ to be a rational elliptic surface with type $IV$ fibres above $t = \zeta$ and $t = \zeta^2$ (the ramification points of $\phi$), and other singular fibres of type $I_2$, $I_1$ and $I_1$. For $X'$ we take a rational surface with an $I_9$ above the point where $X$ has the $I_2$, and three $I_1$'s, two of which above the points where $X$ has $I_1$'s. Now $W$ has 6 singular fibres of type $I_1 \times I_1$, three of type $I_2 \times I_9$ and three of type $I_0 \times I_1$. The $IV$ fibres have become nonsingular in the base change.

From [5] one has that $h^2(X) = h^2(X') = h^2(\mathcal{E}) = 10$ and $h^2(\mathcal{E}') = 34$. From lemma 3.1.1 it follows that $\chi(\widetilde{W}) = 240$. And assuming 3.1.4 one easily computes that $h^2(\widetilde{W}) = 127$. Poincaré duality gives $h^i(\widetilde{W}) = h^{6-i}(\widetilde{W})$. And using that

$$\chi(\widetilde{W}) = \sum_{i=0}^{6}(-1)^i h^i(\widetilde{W})$$

it follows that $h^3(\widetilde{W}) = 16$.

We are also interested in the dimensions of the eigenspaces of the $\sigma$-action on the cohomology. First we compute $\dim H^2(\widetilde{W})_1$. Note that $h^2(\mathcal{E}'/\langle\sigma\rangle) = 10$ and in section 2.3 we showed that $h^2(\mathcal{E}/\langle\sigma\rangle) = 6$. The intersection $r^*(H^2(\mathcal{E})) \cap r'^*(H^2(\mathcal{E}'))$ is generated by the class of a fibre of $\tilde{f}$, and $\sigma$ acts trivial on this class, hence $r^*(H^2(\mathcal{E}))_1 \cap r'^*(H^2(\mathcal{E}'))_1$ is also generated by a fibre class. So we have

$$\begin{aligned}\dim\left(r^*(H^2(\mathcal{E}))_1 + r'^*(H^2(\mathcal{E}'))_1\right) &= \dim H^2(\mathcal{E})_1 + \dim H^2(\mathcal{E}')_1 - 1 \\ &= h^2(\mathcal{E}/\langle\sigma\rangle) + h^2(\mathcal{E}'/\langle\sigma\rangle) - 1 \\ &= 6 + 10 - 1 = 15.\end{aligned}$$

The 84 divisors mentioned in 3.1.4 form 28 orbits of length 3 under the action of $\sigma$. So $\dim H^2(\widetilde{W})_1 = 15 + 28 = 43$.

Since $H^4(\widetilde{W})$ is poincaré-dual to $H^2(\widetilde{W})$ its 1-eigenspace also has dimension 43.

From lemma 3.1.1 is follows that $\chi(\widetilde{W}/\langle\sigma\rangle) = 80$

One has $H^i(\widetilde{W})_1 \cong H^i(\widetilde{W}/\langle\sigma\rangle)$. Combining this with

$$\chi(\widetilde{W}/\langle\sigma\rangle) = \sum_{i=0}^{6}(-1)^i h^i(\widetilde{W}/\langle\sigma\rangle)$$

gives

$$\dim H^3(\widetilde{W})_1 = 2 + 2h^2(\widetilde{W}/\langle\sigma\rangle) - \chi(\widetilde{W}/\langle\sigma\rangle) = 8$$

and

$$\dim H^3(\widetilde{W})_\zeta = \frac{1}{2}(16 - 8) = 4.$$

From lemma 3.1.2 it follows that $h^{3,0}(\widetilde{W}) = 3$ and $h^{3,0}(\widetilde{W}/\langle\sigma\rangle) = 1$ hence $\dim H^{3,0}(\widetilde{W})_\zeta = 1$, and consequently $\dim H^{2,1}(\widetilde{W}/\langle\sigma\rangle)_\zeta = \dim H^{1,2}(\widetilde{W}/\langle\sigma\rangle)_\zeta = \dim H^{0,3}(\widetilde{W}/\langle\sigma\rangle)_\zeta = 1$.

## 3.2 Computation of traces of Frobenius

We are interested in computing $\text{Tr}(\text{Fr}_q^{(i)}|H^3(\widetilde{W}))$. For this we use 2.2.1 again, but now we also need to know $\text{Tr}(\text{Fr}_q^{(i)}|H^2(\widetilde{W}))$ and $\text{Tr}(\text{Fr}_q^{(i)}|H^4(\widetilde{W}))$. The next lemma relates these two traces.

**Lemma 3.2.1**

$$q\,\text{Tr}(\text{Fr}_q^{(i)}(\widetilde{W})) = \text{Tr}(\text{Fr}_q^{(i)}|H^4(\widetilde{W})).$$

*Proof.* Suppose $\mathrm{Fr}_q^{(i)}|H^2(\widetilde{W})(1)$ has eigenvalues $\lambda_j$, $1 \leq j \leq h^2(\widetilde{W})$. We have $|\lambda_j| = 1$. The Poincaré duality pairing

$$H^2(\widetilde{W})(1) \otimes H^4(\widetilde{W})(2) \to H^6(\widetilde{W})(3) = \mathbb{Q}_\ell$$

is $\sigma$-equivariant and Galois equivariant, and $\mathrm{Fr}_q^{(i)}$ acts trivial on $H^6(\widetilde{W})(3)$ (lemma 2.2.1, $(d)$) so the eigenvalues of $\mathrm{Fr}_q^{(i)}|H^4(\widetilde{W})(2)$ are $\lambda_j^{-1}$. Consequently

$$\mathrm{Tr}(\mathrm{Fr}_q^{(i)}|H^2(\widetilde{W})(1)) = \sum_{j=1}^{h^2(\widetilde{W})} \overline{\lambda_j^{-1}} = \overline{\mathrm{Tr}(\mathrm{Fr}_q^{(i)}|H^4(\widetilde{W})(2))} = \mathrm{Tr}(\mathrm{Fr}_q^{(i)}|H^4(\widetilde{W})(2))$$

(the last equality follows from 2.2.1, $(c)$) and the lemma follows. □

For computing the trace on $H^2(\widetilde{W})$ we can use the explicit description given in the previous section. Both Frobenius and $\sigma$ permute the divisors given in 3.1.4, (2), and the trace of $\mathrm{Fr}_q^{(i)}$ is just $q$ times the number of divisors fixed by this map. The traces on $r^*(H^2(\mathcal{E}))$ and $r'^*(H^2(\mathcal{E}'))$ can be computed as explained in section 2.2.

For $q$ big enough one can avoid the tedious computation of the trace on the algebraic part of $H^2(\widetilde{W})$. Assume one knows the trace of the $\mathrm{Fr}_q$-action on the transcedental pieces $H^2(\widetilde{W})^{tr}$ and $H^4(\widetilde{W})^{tr}$ (for example by using lemma 2.2.3). The traces on the algebraic parts satisfy

$$\mathrm{Tr}(\mathrm{Fr}_q^{(i)}|H^2(\widetilde{W})^{alg}) = k\,q\,, \quad \mathrm{Tr}(\mathrm{Fr}_q^{(i)}|H^4(\widetilde{W})^{alg}) = k\,q^2$$

for some integer $k$. The eigenvalues of $\mathrm{Fr}_q^{(i)}|H^3(\widetilde{W})$ have absolute value $q^{3/2}$, so $\left|\mathrm{Tr}(\mathrm{Fr}_q^{(i)}|H^3(\widetilde{W}))\right| \leq h^3(\widetilde{W})q^{3/2}$. Combining this with 2.2.1 one has

$$\big|1 + q^3 + k(q+q^2) + \mathrm{Tr}(\mathrm{Fr}_q^{(i)}|H^2(\widetilde{W})^{tr}) + \mathrm{Tr}(\mathrm{Fr}_q^{(i)}|H^4(\widetilde{W})^{tr})$$
$$-\mathrm{Fix}(\mathrm{Fr}_q^{(i)}|\widetilde{W}(\mathbb{F}_q))\big| \leq h^3(\widetilde{W})q^{3/2}.$$

For $q$ big enough there will be a unique integer $k$ which satisfies this inequality. And knowing $k$, one can compute $\mathrm{Tr}(\mathrm{Fr}_q^{(i)}|H^3(\widetilde{W}))$.

Define $Y = X \times_{\mathbb{P}_t^1} X'$ and let $\widetilde{Y}$ be the blowup of $Y$ at the nodes. Let $\widetilde{Y}_t$ denote the fibre of $\widetilde{Y} \to \mathbb{P}_t^1$ at $t$. The following lemma is analogous to lemma 2.2.4. It can be proved in a similar way.

**Lemma 3.2.2** *Define $K^{(i)}$ as in lemma 2.2.4 and $R^{(i)}$ as the total number points in the fibres of $\widetilde{W}$ above $-\zeta$ and $-\zeta^2$ that are fixed by $\mathrm{Fr}_q^{(i)}$. Then*

$$\mathrm{Fix}(\mathrm{Fr}_q^{(i)}|\widetilde{W}) = 3 \sum_{t \in K^{(i)}} \widetilde{Y}_t(\mathbb{F}_q) + R^{(i)}.$$

Almost all fibres $\widetilde{Y}_t$ are products $X_t \times X'_t$, and we can count points on these fibres by counting points on both factors.

In order to count points on the other fibres we should study the exceptional divisor in the blowup $\widetilde{Y} \to Y$. What happens is that each node of $Y$ gets replaced by its projectivised tangent cone. Geometrically this is a $\mathbb{P}^1 \times \mathbb{P}^1$, but there are two possible Galois actions on such a divisor above an $\mathbb{F}_q$-rational node. Either its ruling is defined over $\mathbb{F}_q$, in which case it has $(q+1)^2$ rational points, or its ruling is not defined over $\mathbb{F}_q$, and it has $q^2 + 1$ rational points.

On an $I_n \times I_m$ fibre there are $nm$ nodes. If both $I_n$ and $I_m$ are split multiplicative fibres then all nodes and all rulings are $\mathbb{F}_q$-rational. If at least one of them (say $I_n$) is non-split and $n$ is even then there are no $\mathbb{F}_q$-rational nodes. If $I_n$ is non-split, n is odd and $I_m$ is split then there are $m$ nodes $\mathbb{F}_q$-rational, and the rulings are not $\mathbb{F}_q$-rational. Finally, if both $I_n$ and $I_m$ are non-split, with $n$ and $m$ odd then there is one $\mathbb{F}_q$-rational node, with $\mathbb{F}_q$-rational ruling.

### 3.3 Equations

In this section we construct explicit equations for $X$ and $X'$ with the required singular fibre configuration.

For $X'$ we can start with a Weierstrass equation for Beauville's elliptic modular surface

$$y^2 = x^3 - \tilde{t}(\frac{1}{3}\tilde{t}^3 + 8)x - \frac{2}{27}\tilde{t}^6 - \frac{8}{3}\tilde{t}^3 - 16 \tag{3}$$

(see [1]) which has singular fibres of type $I_9$ above $\infty$, and of type $I_1$ above $-3$, $-3\zeta$ and $-3\zeta^2$. Later on we will adjust the coordinate $\tilde{t}$.

For $X$ we start with the general Weierstrass equation for a rational elliptic surface with type $IV$-fibres above $\zeta$ and $\zeta^2$:

$$y^2 = x^3 + c_1(t^2 + t + 1)^2 x + (c_2 t^2 + c_3 t + c_4)(t^2 + t + 1)^2. \tag{4}$$

This equation has discriminant $-16(t^2 + t + 1)^4 \Delta(t)$ with

$$\begin{aligned}\Delta(t) &= \left(4 c_1^3 + 27 c_2^2\right) t^4 + \left(8 c_1^3 + 54 c_2 c_3\right) t^3 + \\ &\quad \left(27 c_3^2 + 12 c_1^3 + 54 c_2 c_4\right) t^2 + \left(54 c_3 c_4 + 8 c_1^3\right) t + 4 c_1^3 + 27 c_4^2.\end{aligned}$$

We are going to impose an $I_2$ at $\infty$. For this it is necessary that the $t^4$ and $t^3$ coefficients of $\Delta(t)$ vanish. So we set

$$c_1 = -3c_5^2, \quad c_2 = 2c_5^3, \quad c_3 = 2c_5^3,$$

and $\Delta(t)$ specializes to

$$\left(3456 c_5^6 - 1728 c_4 c_5^3\right) t^2 + \left(3456 c_5^6 - 1728 c_4 c_5^3\right) t + 1728 c_5^6 - 432 c_4^2. \tag{5}$$

There are two $I_1$-fibres above the values of $t$ for which (5) is zero. At these values of $t$ the surface $X'$ will also have to have $I_1$-fibres. In equation (3) two of the three $I_1$-fibres are conjugate in the extension $\mathbb{Q}(\zeta)/\mathbb{Q}$. So this leads us to require that the discriminant of (5) with respect to $t$ is $-3$ times a square. This discriminant is

$$-2^{12}3^6 \left(c_5{}^3 + c_4\right) c_5{}^3 \left(2\,c_5{}^3 - c_4\right)^2,$$

and we make it $-3$ times a square by setting $c_4 = 3c_6^2 c_5 - c_5^3$. Now (5) becomes

$$81c_5^2(c_6 - c_5)(c_6 + c_5)(4c_5^2 t^2 + 4c_5^2 t + c_5^2 + 3c_6^2). \tag{6}$$

We want two $I_1$-fibres of $X'$ to lie above the zeroes of (6), keeping the $I_9$ at $\infty$. This can be done with a lineair change of the coordinate $\tilde{t}$ in (3). Set

$$\tilde{t} = \frac{6\,c_5\,t + 3\,c_5 + 3\,c_6}{2\,c_6}.$$

Now we have a family of elliptic surfaces $X'$ defined by (3), and $X$ defined by (4) with all the substitutions made. For $c_5$ and $c_6$ general enough they will satisfy the required singular fibre condition.

For future $L$-function calculations it is convenient to construct a Galois representation with only few ramified primes. For this we want our threefold $\widetilde{W}$ to have only a few primes of bad reduction. The primes of bad reduction are precisely the primes modulo which different singular fibres coincide. So we would like to specialize $c_5$ and $c_6$ such that this set of primes is small.

The $t$-values at which there is a singular fibre are $t = \infty$ and $t$ is a zero of

$$(3\,c_6 + c_5 + 2\,c_5 t)\left(3\,c_6{}^2 + 4\,c_5{}^2 t^2 + 4\,c_5{}^2 t + c_5{}^2\right)\left(t^2 + t + 1\right). \tag{7}$$

The discriminant of (7) with respect to $t$ is

$$2^{12} 3^{10} c_5{}^6 c_6{}^6 \left(3\,c_6{}^2 + c_5{}^2\right)^2 (c_6 - c_5)^4 (c_6 + c_5)^4$$

so the primes which, after specialisation of $c_5$ and $c_6$, occur in this number are the primes of bad reduction. A good choice for $c_5$ and $c_6$ is $c_5 = 3$, $c_6 = 1$. Now only the primes 2 and 3 are bad. The Weierstrass equations in this case are

$$\begin{array}{rcl}
X & : y^2 & = x^3 - 27\,(t^2 + t + 1)^2 x + 18\,(3\,t^2 + 3\,t - 1)(t^2 + t + 1)^2, \\
X' & : y^2 & = x^3 - 3\,(3\,t + 2)(243\,t^3 + 486\,t^2 + 324\,t + 80)\,x \\
& & \quad -39366\,t^6 - 157464\,t^5 - 262440\,t^4 - 235224\,t^3 \\
& & \quad -120528\,t^2 - 33696\,t - 4048.
\end{array} \tag{8}$$

### 3.4 Computation of characteristic polynomials

Define

$$\rho\ :\ G_\mathbb{Q}\ \longrightarrow\ \mathrm{Aut}(H^3(\widetilde{W})_\zeta).$$

Let $p$ be a prime where $\widetilde{W}$ has reduction, and let $a_p = \mathrm{Tr}(\rho(\mathrm{Fr}_p))$ and $b_p = \mathrm{Tr}(\rho(\mathrm{Fr}_{p^2}))$. Then the characteristic polynomial of $\mathrm{Fr}_p$ has the form

$$T^4 + d_3 T^3 + d_2 T^2 + d_1 T + d_0$$

with
$$d_3 = -a_p \quad \text{and} \quad d_2 = \frac{1}{2}(a_p^2 - b_p).$$

With lemma 2.1.1 we can compute the other coefficients. We find
$$\xi = d_2/\bar{d}_2, \quad d_1 = \xi p^3 d_3 \quad \text{and} \quad d_0 = \xi p^6.$$

The first few characteristic polynomials are listed in the next table.

| 5  | $T^4 + (10 + 13\zeta)T^3 - 5\zeta^2 T^2 + 5^3(13 + 10\zeta)T + 5^6\zeta$ |
|----|---|
| 7  | $T^4 - (7 + 4\zeta)T^3 - 189\zeta T^2 + 7^3(7 + 3\zeta)T + 7^6\zeta^2$ |
| 11 | $T^4 - (21 + 19\zeta)T^3 + 517\zeta T^2 + 11^3(21 + 2\zeta)T + 11^6\zeta^2$ |
| 13 | $T^4 + (77 + 70\zeta)T^3 - 1742\zeta^2 T^2 + 13^3(70 + 77\zeta)T + 13^6\zeta$ |
| 17 | $T^4 + (24 - 63\zeta)T^3 - 1802\,T^2 + 17^3(87 + 63\zeta)T + 17^6$ |
| 19 | $T^4 + (73 + 81\zeta)T^3 - 4275\,T^2 - 19^3(8 + 81\zeta)T + 19^6$ |
| 23 | $T^4 + (33 + 129\zeta)T^3 + 14536\zeta^2 T^2 + 23^3(129 + 33\zeta)T + 23^6\zeta$ |
| 29 | $T^4 - (186 + 100\zeta)T^3 + 16936\zeta T^2 + 29^3(186 + 86\zeta)T + 29^6\zeta^2$ |

The determinant character of $\rho$ is the cyclotomic character times a Dirichlet character with values in $\mathbb{Q}(\zeta)$. The only primes that can ramify in the Dirichlet character are 2 and 3, hence its conductor is a divisor of 72. From the table above it follows that it is in fact a character of order 3 and conductor 9. Denote the representation $\rho$ twisted with this Dirichlet character by $\tilde{\rho}$. The determinant of $\tilde{\rho}$ is the cyclotomic character.

## 4 Computation of the $L$-function

We will compute the $L$-function of $\tilde{\rho}$, and not of $\rho$. The reason for this is that we expect $\tilde{\rho}$ to have a smaller conductor than $\rho$ and that makes the computation easier. The $L$-function is an Euler product of local $L$-factors
$$L(\tilde{\rho}, s) = \prod_{p \text{ prime}} L_p(\tilde{\rho}, s).$$

It converges absolutely if $\text{Re}(s) > \frac{5}{2}$. For unramified primes the local factors are defined by
$$L_p(\tilde{\rho}, s) = \frac{1}{p^4 P_p(p^{-s})}$$

with $P_p(T)$ the characteristic polynomial of $\tilde{\rho}(\text{Fr}_p)$.

For ramified primes $p$ we look at the subspace where the inertia group $I_{\mathfrak{p}}$ of a prime $\mathfrak{p}$ of $\bar{\mathbb{Q}}$ above $p$ acts trivial. Let $P_p(T)$ denote the characteristic polynomial of $\text{Fr}_p$ acting on the inertia invariants. Then
$$L_p(\tilde{\rho}, s) = \frac{1}{p^{\deg P_p(T)} P_p(p^{-s})}.$$

Expanding the Euler product we get a Dirichlet series

$$L(\tilde{\rho}, s) = \sum_{n=1}^{\infty} a_n n^{-s}.$$

To compute $a_n$ for all $n$ less than some bound $B$ we need to compute $a_p = \text{Tr}(\tilde{\rho}(\text{Fr}_p))$ for $p < B$ and $b_p = \text{Tr}(\tilde{\rho}(\text{Fr}_{p^2}))$ for $p < \sqrt{B}$.

Let $\tilde{\rho}^*$ be the contragredient representation of $\tilde{\rho}$. Its Tate twist $\tilde{\rho}^*(-3)$ has eigenvalues of Frobenius of the same absolute value as $\tilde{\rho}$. From lemma 2.1.1 it follows that its $L$-function is

$$L^*(s) \;=\; L(\tilde{\rho}^*(-3), s) \;=\; \sum_{n=1}^{\infty} \bar{a}_n n^{-s}.$$

Define the completed $L$-function $\Lambda(\tilde{\rho}, s)$ by

$$\Lambda(s) = \Lambda(\tilde{\rho}, s) = N^{s/2}(4\pi^2)^{-s}\Gamma(s)\Gamma(s-1)L(\tilde{\rho}, s), \tag{9}$$

with $N$ the conductor of $\tilde{\rho}$, and define $\Lambda^*(s) = \Lambda(\tilde{\rho}^*(-3), s)$. It is conjectured that $\Lambda(s)$ can be analytically continued to the whole complex plane, and that it satisfies the functional equation

$$\Lambda(s) \;=\; w\,\Lambda^*(4-s) \tag{10}$$

for some $w \in \mathbb{C}$, $|w| = 1$.

There is a well known trick for computing $L(2)$ if one assumes the functional equation (cf.[2], but note that in this paper, page 119, one should take $c = \sqrt{N}/4\pi^2$ instead of the erroneous $c = \sqrt{N}/16\pi^2$). The idea is to integrate

$$\frac{\Lambda(2+s)\,t^{-s}}{2\pi i\, s}$$

along a path in $\mathbb{C}$ around 0 for some number $t$. By Cauchy's theorem this integral is $\Lambda(2)$. One can stretch the path more and more in such a way that in the limit it besomes a union of two lines, one going from $r - i\infty$ to $r + i\infty$ and the other going from $-r + i\infty$ to $-r - i\infty$, for some real $r > \frac{1}{2}$. This can be done because the $\Gamma$-factors in (9) make sure that the integral over pieces of paths with very big or very small imaginary part is very small. So one has

$$\Lambda(2) = \frac{1}{2\pi i}\int_{r-i\infty}^{r+i\infty}\frac{\Lambda(2+s)t^{-s}}{s}ds + \frac{1}{2\pi i}\int_{-r+i\infty}^{-r-i\infty}\frac{\Lambda(2+s)t^{-s}}{s}ds. \tag{11}$$

Define

$$F(x) = \frac{1}{2\pi i}\int_{r-i\infty}^{r+i\infty}\Gamma(s+2)\Gamma(s)x^{-s}ds.$$

See [2] for a discussion on how to compute values of $F$ efficiently.

The first integral of (11) is equal to

$$\frac{1}{2\pi i} \int_{r-i\infty}^{r+i\infty} \Gamma(s+2) \frac{\Gamma(s+1)}{s} \left(\frac{4\pi^2}{\sqrt{N}}\right)^{-s-2} t^{-s} \sum_{n=1}^{\infty} \frac{a_n}{n^{s+2}} ds =$$

$$\frac{N}{16\pi^4} \sum_{n=1}^{\infty} \frac{a_n}{n^2} \frac{1}{2\pi i} \int_{r-i\infty}^{r+i\infty} \Gamma(s+2)\Gamma(s) \left(\frac{4nt\pi^2}{\sqrt{N}}\right)^{-s} ds =$$

$$\frac{N}{16\pi^4} \sum_{n=1}^{\infty} \frac{a_n}{n^2} F\left(\frac{4nt\pi^2}{\sqrt{N}}\right)$$

Using the functional equation (10) one can rewrite the second integral of (11) as

$$\frac{w}{2\pi i} \int_{r-i\infty}^{r+i\infty} \frac{\Lambda^*(2+s)t^s}{s} = \frac{wN}{16\pi^4} \sum_{n=1}^{\infty} \frac{\bar{a}_n}{n^2} F\left(\frac{4n\pi^2}{t\sqrt{N}}\right).$$

Consequently one has

$$L(2) = \sum_{n=1}^{\infty} \frac{a_n}{n^2} F\left(\frac{4nt\pi^2}{\sqrt{N}}\right) + w \sum_{n=1}^{\infty} \frac{\bar{a}_n}{n^2} F\left(\frac{4n\pi^2}{t\sqrt{N}}\right). \tag{12}$$

We will use (12) to compute $L(2)$, or at least obtain a likely candidate for $L(2)$. But there are some problems. We are only able to compute Euler factors of $L$ at primes where $\widetilde{W}$ has good reduction. And we don't know $N$ and $w$.

We proceed as follows. We guess the Euler factors at bad primes. Assuming conjecture $C_6$ in [8] it follows there are only finitely many possibilities. Moreover, we guess the conductor. The only primes that can occur in the conductor are the primes of bad reduction. Note the left hand side of (12) is independent of $t$, so the right hand side is independent too. We compute the right hand side for several values of $t$, and from the fact that it is constant we get a system of equations in $w$ (with as many equations as we like). If we made wrong guesses for bad Euler factors and conductor the system will probably have no solution, and we make new guesses. We continue this way until there appears to be a $w$ which satisfies all equations, and moreover has absolute value close to 1.

For $\widetilde{W}$ given by equations (8) we computed $a_n$ for $n < 80000$. For Euler factors at 2 and 3 both equal to 1, and conductor $2^9 3^9$ the above procedure was performed for $t = 1, 1.2, 1.3, 1.4$ and $1.5$. The right hand side of (12) remained constant upto 20 decimal places, and the value of $w$ that came out was upto 20 decimal places equal to a sixth root of unity $-\zeta$. Replacing $w$ with $-\zeta$ the value of (12) remained constant upto 24 decimal places. With this choice we found

$$L(2) = 0.419940577402445298239298 4 + 0.242452805406948667234606 3\, i.$$

Jasper Scholten
Vakgroep Wiskunde
Rijksuniversiteit Groningen
postbus 800
9700 AV Groningen
the Netherlands
e-mail: jasper@math.rug.nl